\magnification=1200

\font\sets=msbm10.
\font\stampatello=cmcsc10.

\def\square{\hbox{\vrule\vbox{\hrule\phantom{s}\hrule}\vrule}}
\def\defineq{\buildrel{def}\over{=}}
\def\N{\hbox{\sets N}}
\def\R{\hbox{\sets R}}
\def\dIk{\Delta_{_{{\cal J}_k}}}

\par
\centerline{\bf A note on the exponential sums of the localized
divisor functions}

\medskip
\par
\centerline{\stampatello giovanni coppola and maurizio laporta}
\footnote{}{\par \noindent MSC(2010): 11L07,11N99. Keywords: {\it exponential sum, divisor function}}
\medskip
\par
\noindent
For a fixed integer $k\ge 2$, let us consider $k-1$ intervals $I_j\subseteq \N$ ($j=1,\ldots, k-1$) and define
${\cal J}_k\defineq I_1\times \cdots \times I_{k-1}$. Then, for all  $n\in\N$ let us set
$$
\dIk(n)\defineq\#\{(n_1,\ldots,n_{k-1})\in {\cal J}_k:\ n_1\cdots n_{k-1}|n\}.
$$
We say that $\dIk$ is the {\it divisor function localized} on ${\cal J}_k$. 
By taking ${\cal J}_k=\N^{k-1}$ we recover the standard divisor function $d_k(n)$. Moreover, if all the $I_j$ are intervals of logarithmic length $1$, then 
$\max_{{\cal J}_k}\dIk(n)$ is the {\it concentration} function introduced by Hooley [Ho] (see also [HT]).
Given 
$\alpha\in[0,1)$ and $N\in\N$, the exponential sum associated to $\dIk$ over the integers of $(N,2N]$ is
$$
S_k(\alpha,N)\defineq \sum_{n\sim N}\dIk(n)e(n\alpha), 
$$
\par
\noindent
where  $n\sim N$ means that $n\in(N,2N]\cap\N$ and we write $e(\beta)$ for $e^{2\pi i \beta}$. Throughout, the Vinogradov notation $\ll_{k,\varepsilon}$ is synonymous of Landau's $O_{k,\varepsilon}$ ($\varepsilon>0$ is arbitrarily small and may change at each occurrence). 
\smallskip
\par
\noindent
{\stampatello theorem.} {\it For all relatively prime integers $a,q$  with $q>1$, one has, uniformly in \thinspace $a$, 
$$
S_k\Big({a\over q},N\Big)\ll_{k,\varepsilon}(Nq)^{\varepsilon}\Big({N\over q}+q+N^{1-1/k}\Big).
$$
\par
\noindent
Furthermore, this same upper bound applies to $S_k(\alpha,N)$ provided that $\alpha\in (0,1)$ satisfies
$|\alpha-a/q|\le 1/q^2$.}
\smallskip
\par
\noindent
{\stampatello proof.} First, let us assume that 
$\{n\in(N,2N]\cap\N: \dIk(n)\not=0\}\not=\emptyset$ for otherwise the inequality is trivial.
Then, we denote ${\cal I}_k\defineq {\cal J}_k\times I_{k}$ with
$$
I_{k}\defineq\{m\in\N: n_1\cdots n_{k-1}m\sim N\ \hbox{for some}\ (n_1,\ldots,n_{k-1})\in {\cal J}_k\}.
$$
Moreover, let us
write $S_k(\alpha,N)$ as a multiple exponential sum,
$$
S_k(\alpha,N)
=\sum_{{\vec{n}_k\in{\cal I}_k}\atop {{\bf n}_k\sim N}}e({\bf n}_k\alpha),
$$
where we set \thinspace $\vec{n}_k\defineq (n_1,\ldots,n_{k})$ \thinspace and \thinspace ${\bf n}_k\defineq n_1\cdots n_k$, for brevity. Since it is plain that 
${\bf n}_k>N$ implies that 
$n_j>N_k\defineq[N^{1/k}]$ for some $j\in\{1,\ldots,k\}$ (hereafter $[x]$ is the integer 
				% PAGE 2 
part of $x\in\R$), we write
$$
\eqalign{
S_k(\alpha,N)=&\sum_{{{\vec{n}_k\in{\cal I}_k}\atop {{\bf n}_k\sim N}}\atop{n_1> N_k}}e({\bf n}_k\alpha)
 +\sum_{{{\vec{n}_k\in{\cal I}_k}\atop {{\bf n}_k\sim N}}\atop{n_1\le N_k}}e({\bf n}_k\alpha)\cr
=&\sum_{{{\vec{n}_k\in{\cal I}_k}\atop {{\bf n}_k\sim N}}\atop{n_1> N_k}}e({\bf n}_k\alpha)
+\sum_{{{\vec{n}_k\in{\cal I}_k}\atop {{\bf n}_k\sim N}}\atop{n_1\le N_k<n_2}}e({\bf n}_k\alpha)
+\sum_{{{\vec{n}_k\in{\cal I}_k}\atop {{\bf n}_k\sim N}}\atop{n_1,n_2\le N_k}}e({\bf n}_k\alpha)\cr
=&\sum_{{{\vec{n}_k\in{\cal I}_k}\atop {{\bf n}_k\sim N}}\atop{n_1> N_k}}e({\bf n}_k\alpha)
+\sum_{{{\vec{n}_k\in{\cal I}_k}\atop {{\bf n}_k\sim N}}\atop{n_1\le N_k<n_2}}e({\bf n}_k\alpha)+\sum_{{{\vec{n}_k\in{\cal I}_k}\atop {{\bf n}_k\sim N}}\atop{n_1,n_2\le N_k<n_3}}e({\bf n}_k\alpha)+ \sum_{{{\vec{n}_k\in{\cal I}_k}\atop {{\bf n}_k\sim N}}\atop{n_1,n_2,n_3\le N_k}}e({\bf n}_k\alpha)\cr
=&\ldots\ldots\cr
=&\sum_{{{\vec{n}_k\in{\cal I}_k}\atop {{\bf n}_k\sim N}}\atop{n_1> N_k}}
e({\bf n}_k\alpha)+\sum_{j=2}^{k}\sum_{{{\vec{n}_k\in{\cal I}_k}\atop {{\bf n}_k\sim N}}\atop{n_1,\ldots,n_{j-1}\le N_k<n_j}}
e({\bf n}_k\alpha).\cr}
$$
\par
\noindent
Of course, it is tacitly understood that if the constraints $n_1,\ldots,n_{j-1}\le N_k$, $N_k<n_j$ are incompatible respectively with 
$(n_1,\ldots,n_{j-1})\in I_1\times \cdots \times I_{j-1}$, $n_j\in I_{j}$, then the sum under any of such conditions is meant to be zero.
Thus, we have
$$
\eqalign{
|S_k(\alpha,N)|\leq&\sum_{j=1}^{k}\sum_{{\vec{n}_k^{(j)}\in {\cal I}_k^{(j)}}\atop {{\bf n}_k/n_j<2N/N_k}}
\Big|\sum_{{n_j\in I_j}\atop {    
                              {{\bf n}_k\sim N}\atop {n_j>N_k}
                              }
           }e({\bf n}_k\alpha)\Big|, \cr}
$$
\par
\noindent
where ${\cal I}_k^{(j)}\defineq I_1\times\ldots \times I_{j-1}\times I_{j+1}\times \ldots \times I_{k}$
and analogously $\vec{n}_k^{(j)}$ is the $k-1$ dimensional vector obtained by removing the $j$-th entry from $\vec{n}_k$. We apply the well-known inequality [D, Ch.25] 
$$
\Big|\sum_{{n_j\in I_j}\atop {    
                              {{\bf n}_k\sim N}\atop {n_j>N_k}
                              }
           }e({\bf n}_k\alpha)\Big|
=\Big|\sum_{{n_j\in I_j,n_j>N_k}\atop {n_j\sim {Nn_j\over{\bf n}_k}}}e\Big(\alpha{{\bf n}_k\over n_j}n_j\Big)\Big|
 \le \min\Big({Nn_j\over {\bf n}_k},{1\over\Vert\alpha{\bf n}_k/n_j\Vert}\Big),
$$
\par
\noindent
where $\Vert x\Vert$ denotes the distance of $x\in\R$ from the integers.
Thus, writing $t={\bf n}_k/n_j$, 
$$
|S_k(\alpha,N)|
\le k\sum_{t<2N/N_k}d_{k-1}(t)\min\Big({N\over t},{1\over\Vert t\alpha\Vert}\Big). 
$$
\par
\noindent
Note that $d_{1}(t)=1$, while for $k\ge 2$ recall that $d_{k}(t)\ll_{k,\varepsilon}t^{\varepsilon}, \forall \varepsilon>0$. Hence, by taking 
				% PAGE 3
$\alpha=a/q$ and denoting with $\overline{a}$ the inverse of $a$ $(\bmod\, q)$, we get
$$
\eqalign{
S_k\Big({a\over q},N\Big)
\ll_{k,\varepsilon}&
N^{\varepsilon}\Big({N\over q}\sum_{t'<{{2N}\over {qN_k}}}{1\over {t'}}+\sum_{1\le r\le {q\over 2}}{q\over r}\sum_{{t<2N/N_k}\atop {t\equiv \pm r\overline{a}\, (\!\! \bmod q)}}1\Big)\cr
\ll_{k,\varepsilon}
&N^{\varepsilon}\Big({N\over q}+\sum_{1\le r\le {q\over 2}}{q\over r}\Big({N\over {qN_k}}+1\Big)\Big) 
\cr
\ll_{k,\varepsilon}&(Nq)^{\varepsilon}\Big({N\over q}+{N\over {N_k}}+q\Big)\cr
\ll_{k,\varepsilon}&(Nq)^{\varepsilon}\Big({N\over q}+q+N^{1-1/k}\Big).}
$$
\par
\noindent
Once $\alpha\in (0,1)$ is such that 
$|\alpha-a/q|\le 1/q^2$, 
the bound for $S_k(\alpha,N)$ follows by the same calculations to prove $(3)$ in [D,Ch.25]. The Theorem is completely proved.
\hfill \square 

\medskip

\par
\noindent {\bf Remark 1.} It transpires that the upper bound of $S_k(\alpha,N)$ does not depend 
on the {\it localization} of the divisors of $n\in(N,2N]$. In particular, it holds also for the exponential sum associated to the 
Hooley's function. 
\par
For the divisor function $d_k$ we explicitly state the following. 
\smallskip
\par
\noindent
{\stampatello corollary.} {\it For all relatively prime integers $a,q$ with $q>1$ we have, uniformly for $\alpha\in[a/q-1/q^2,a/q+1/q^2]$ and \thinspace $a$,} 
$$
\sum_{n\sim N}d_k(n)e(n\alpha)\ll_{k,\varepsilon}(Nq)^{\varepsilon}\Big({N\over q}+q+N^{1-1/k}\Big). 
$$
\par
\noindent {\bf Remark 2.} Clearly, such an inequality gives some improvement on the trivial bound
$N^{1+\varepsilon}$ when $N^{1/k}\ll q\ll N^{1-1/k}$, which indeed yields the estimate $N^{1-1/k+\varepsilon}$. It is under these conditions that such inequalities are most commonly applied. 
\par
We are going to use them in our treatment of mean-squares of $d_k$  in short intervals (see [CL1] and [CL2]). In particular,  
the present investigation on $\dIk$ has been motivated by
our $k-$folding method [CL1], where we {\it localize} the divisor function in suitable boxes. 

\medskip

\par
\noindent {\bf Acknowledgment.} The authors wish to thank A. Ivi\'c for useful hints. 

\medskip

\par
\centerline{\bf Bibliography}
\smallskip
\item{\bf [CL1]} Coppola, G. and Laporta, M. \thinspace - \thinspace {\sl Generations of correlation averages} \thinspace - \thinspace Journal of Numbers Volume 2014 (2014), Article ID 140840, 1-13 (available at
\par 
http://dx.doi.org/10.1155/2014/140840) 
				% PAGE 4 
\item{\bf [CL2]} Coppola, G. and Laporta, M. \thinspace - \thinspace {\sl Sieve functions in arithmetic bands} \thinspace - \thinspace preprint at 
http://arxiv.org/abs/1503.07502v2 
\item{\bf [D]} Davenport, H. \thinspace - \thinspace {\sl Multiplicative Number Theory} \thinspace - \thinspace Third Edition, GTM 74, Springer, New York, 2000 
\item{\bf [HT]} Hall, R.R. and Tenenbaum, G. \thinspace - \thinspace  {\sl Divisors} \thinspace - \thinspace Cambridge Tract 90, Cambridge University Press, Cambridge, (1988)
\item{\bf [Ho]} Hooley, C.\thinspace - \thinspace  {\sl On a new technique and its applications to the theory
  of numbers} \thinspace - \thinspace Proc. London Math. Soc. (3) {\bf 38} (1979), 115-151

\medskip

\par
\leftline{\tt Giovanni Coppola \hfill Maurizio Laporta}
\leftline{\tt Home \negthinspace address:  \hfill Universit\`a degli Studi di Napoli}
\leftline{\tt V.Partenio,12 \hfill Dip.di Matematica e Applicazioni}
\leftline{\tt 83100, Avellino (ITALY) \hfill Complesso di Monte S.Angelo,}
\leftline{\tt www.giovannicoppola.name \hfill Via Cinthia - 80126, Napoli (ITALY)}
\leftline{\tt giovanni.coppola@unina.it \hfill mlaporta@unina.it}

\bye